\documentclass{amsart} 
\usepackage{latexsym,amssymb,graphics,txfonts}
\usepackage[pdftex]{graphicx}
\newtheorem{thm}{Theorem}[section]

\newtheorem{prop}[thm]{Proposition}

\newcommand{\BC}{{\mathbb{C}}}
\newcommand{\BZ}{{\mathbb{Z}}}
\newcommand{\BQ}{{\mathbb{Q}}}
\newcommand{\BR}{{\mathbb{R}}}
\newcommand{\rp}{{\BR P(2)}}

\newcommand{\A}{{\BR A}}
\newcommand{\CA}{{\BC A}}
\newcommand{\cp}{{\BC P(2)}}
\newcommand{\ra}{{\rightarrow}}
\newcommand{\Bc}{{\mathcal{C}}}
\newcommand{\T}{{\mathcal{T}}}
\newcommand{\D}{{\mathcal{D}}}
\newcommand{\I}{{\mathcal{I}}}

\newcommand{\be}{\begin{equation}}
\newcommand{\ee}{\end{equation}}
\DeclareMathOperator{\lk}{lk}
\DeclareMathOperator{\rad}{rad}
\DeclareMathOperator{\Arf}{Arf}
\DeclareMathOperator{\td}{td}

\begin{document}

\title{Arf invariants of Real algebraic curves }

\dedicatory{Dedicated to the memory of my adviser, \\ P. Emery Thomas}

\author{ Patrick M. Gilmer}
\address{Department of Mathematics\\
Louisiana State University\\
Baton Rouge, LA 70803\\
USA}
\email{gilmer@math.lsu.edu}
\thanks{This research was partially supported by NSF-DMS-0203486. Some of this research was conducted while participating in the MSRI Program: Topological Aspects of  Real Algebraic Geometry, Spring 2004}
\urladdr{www.math.lsu.edu/\textasciitilde gilmer/}

\date{November 7,  2005}

\begin{abstract} We give some new congruences for singular real algebraic curves which generalize Fiedler's congruence
for nonsingular curves.\end{abstract}

\maketitle

\section{Introduction}

Let $\A$ be an irreducible real algebraic curve of degree $2k$ with only real nodal singularities. $\A$ consists of the image of a number of immersed circles. $\A$ is called an M-curve if the number of immersed circles plus the number of double points is 
$1 + \binom{2k-1}{2}$. Let $\CA$ be the complex curve in $\cp$ given by the same polynomial as $\A.$ Thus $\A = \CA \cap \rp.$ $\A$ is a M-curve precisely when
$\CA \setminus \A$ consists of two punctured spheres which are interchanged by complex conjugation.  Arbitrarily choose one of these components, say $\CA^+$. The complex structure on $\CA^+$ induces an orientation on $\CA^+$, and thus on each immersed circle of  $\A$.
Of course if we choose the other component,  we would get the opposite orientation on each immersed circle of $\A.$  An orientation on each of the components up to reversing all the orientations simultaneously is called a semi-orientation. Thus each M-curve has receives a semi-orientation, called the complex orientation \cite{R}.

 An oval is a 2-sided simple closed curve in $\rp$, the real projective plane.
The inside of an oval is the component of its complement which is a disk. The outside of an oval is the component of its complement which is a Mobius band. Suppose $\Bc$ is the disjoint collection of oriented ovals.   We call such a collection a simple curve.  An oval of $\Bc$ is called even  if it lies inside and even number of other ovals of $\Bc$.  An oval of $\Bc$ is called odd if it lies inside an odd number of other ovals of $\Bc$. Let $p(\Bc)$ denote the number of even ovals in $\Bc$, and  $n(\Bc)$ denote the number of odd ovals in $\Bc$. If one of ovals of $\Bc$ lies inside another oval of $\Bc$ or vise versa, we say they are linked. We say $\Bc$ is odd if each oval is linked with an odd number of other ovals. It follows that an odd curve must have an even number of components.  We say $\Bc$ is even if each oval is linked with an even number of other ovals and the total number of ovals is odd.

Let $\Pi^+(\Bc)$ denote the number of pairs of linked ovals for which the orientations on the curves extend to an orientation of the intervening annulus.   Let $\Pi^-(\Bc)$ denote the number of pairs of linked ovals for which the orientations on the curves do not extend to an orientation of the intervening annulus.

We will usually write simply $n$, $p$, $\Pi^\pm$ without $\mathcal C$, except in instances where it might be unclear to which curve these numerical characteristics refer.

By a curve,  we will mean a collection of immersed oriented oriented curves in $\rp$
with only transverse double point intersections. 
We will say two curves $A$ and $A'$ are weakly equivalent if they can be connected a sequence of ambient  isotopies in $\rp$, and the local moves and their inverses: balanced type I moves (see Figure \ref{bal} below), safe type II \footnote{This is a Reidemeister type II move on diagrams without over and under crossings, but where the two strands have opposite orientations}, type III \footnote{ This is a Reidemeister type III move on diagrams without over and under crossings}, and  ``empty" figure eight death. 

\begin{figure}[ht]
\includegraphics[width=1.5in]{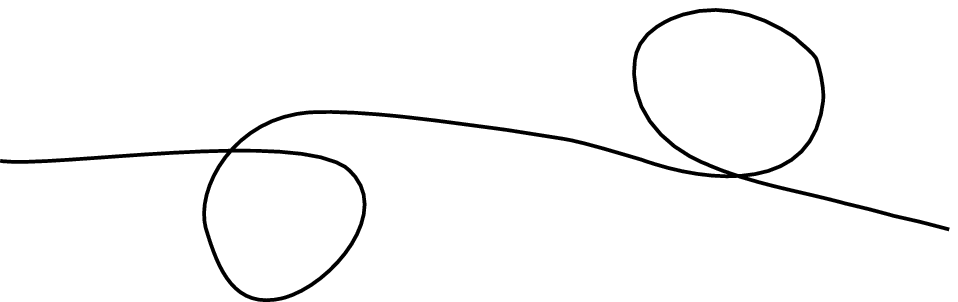}
\caption{A balanced type one move replaces an arc with two curls like that above by an arc without any double points.} \label{bal}
\end{figure}

\begin{thm}\label{simple.thm} Let $\A$ be a nodal M-curve of degree $2k$. 

If $k$ is even, and $\A$ is weakly equivalent to an odd simple curve $\Bc$, then
\[ 
\Pi^+(\Bc) -\Pi^-(\Bc)- p(\Bc) \equiv  
\begin{cases} \frac {k^2} {2}  \pmod{8} & \\ \text{ or } &\\  \frac {k^2} {2} -2  \pmod{8}& \end{cases}.
 \]

If $k$ is odd, and $\A$ is weakly equivalent to an even simple curve $\Bc$, then
\[ 
\Pi^+(\Bc) -\Pi^-(\Bc)- n(\Bc) \equiv  
 \frac {k^2-1 } {2} 
\pmod{8}. \]
\end{thm}

\begin{figure}[ht]
\includegraphics[width=2in]{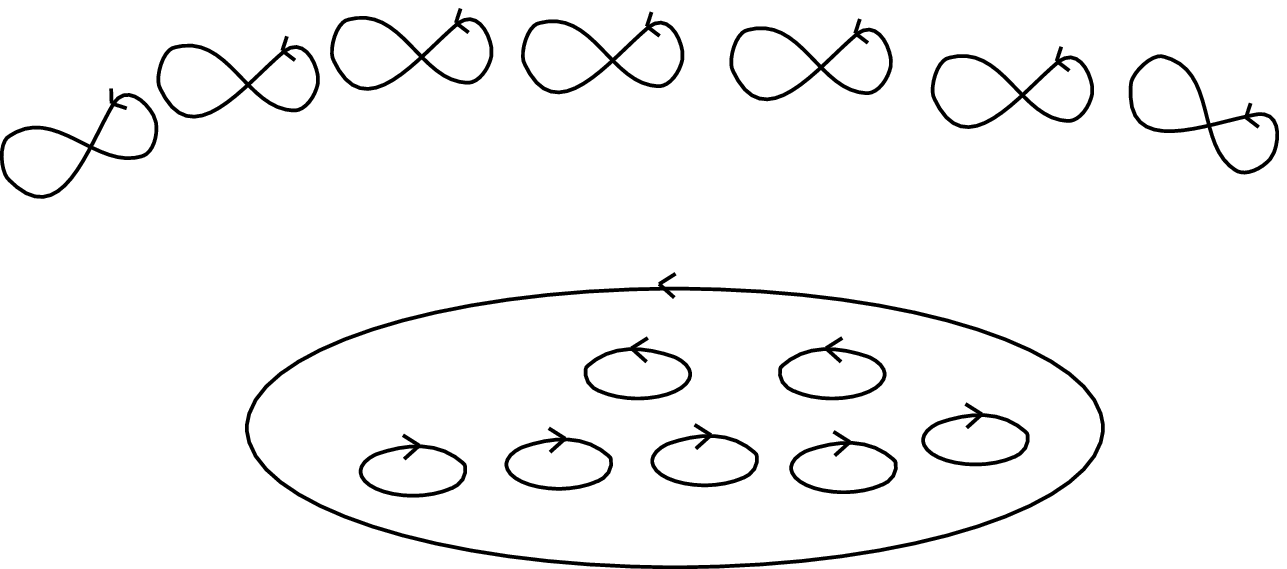}
\caption{A scheme for a degree 8 nodal M-curve with complex orientation prohibited by Theorem \ref{simple.thm}} \label{nodal}
\end{figure}

 In the case $\A= \Bc$,  this becomes a congruence due to Fiedler \cite{F1}. This is not obvious. We prove this explicitly in Theorem \ref{fied} . Fiedler \cite{F2} has also given some generalizations of his congruence to singular curves. The scheme  in Figure \ref{nodal} below is prohibited by Theorem \ref{simple.thm}  but not by the results of \cite{F2}. Nor is it prohibited by any of the other known general restrictions on nodal curves: the Kharlamov-Viro congruences \cite{KV} (as correctly stated in \cite{ OV}), the Viro inequalities \cite{V2, Fi, G5}, and also the extremal properties of the Viro inequalities in \cite{G5}.  If we perturb the figure eights  in  this scheme into  pairs of ovals we obtain the scheme for an nonsingular M-curve $< 14 \coprod 1< 7>>$ which can be realized by a real algebraic curve \cite{V1}.  
 
 Theorem  \ref{simple.thm} is a corollary of Theorem \ref{main} below.  This theorem is actually a simplified version of the more general
 \cite[Theorem 3.3]{G5}.    In order to prove (or even state) Theorem \ref{main}, we must discuss the  Arf invariant of links which we do in sections two and three. These sections are, mainly, a review of parts of \cite{G2,G3}.

It seems likely that any scheme for a curve  with complex orientation which can be prohibited by \cite{F2} can also be prohibited by  our Theorem \ref{simple.thm}. However our theorem applies to hypothetical curves  with a complex orientation. Fiedler's  result concerns hypothetical singular curves which, by hypothesis, are related to actual nonsingular real algebraic curves by desingularization. Thus information about the complex orientations of the hypothetical curves is contained only implicitly in the relation to the actual nonsingular real algebraic curve.
For this reason, it seems difficult to derive \cite{F2} as a corollary 
of Theorem  \ref{main} in a way similar to the proof of Theorem \ref{fied}.  However both obstructions can be interpeted as deriving from the calculation of the Brown invariant of the Gillou-Marin form on  characteristic surfaces. Moreover the two surfaces are closely related.  As evidence for the above suggestion, we will 
show, in the last section,  that the two explicit examples given by Fiedler of curves prohibited by his congruence for singular curves are also prohibited by Theorem \ref{main}.

\section{Arf invariants of links in $S^3$}

An oriented link   in $S^3$ is called proper if for each component the sum of the linking 
numbers with all the other components is even. i.e. $L = \sqcup_i K_i$ is proper if and only if 
$\lk(K_i, L-K_i) \equiv  0 \pmod{2}$ for all $i.$ We use $\lk$ to denote the $\BZ$-valued linking number of oriented links.

Robertello defined the Arf invariants of proper links and gave several equivalent definitions. One involved the Seifert pairing on an orientable spanning surface. 
We generalized this definition so that it applies to non-orientable spanning surfaces
\cite{G1,G2}. This definition for the Arf invariant is analogous to  the Gordon-Litherland 
\cite{GL} definition of the signature of a knot. There is a version for unoriented links
but the version for oriented links is more useful in this paper.

Let $V$ be a $\BZ/2\BZ$ vector space. equipped with a  bilinear symmetric form
\[ \cdot :V \times V  \rightarrow \BZ/2\BZ.\]
A function $q: V \rightarrow   \BZ/4\BZ$ is called a quadratic refinement of $\cdot$, if
\[q(x+y)-q(x) -q(y) = 2 \  x \cdot  y \] 
holds  for all $x, y \in V$.
Here  $2$ denotes the non-trivial group homomorphism $\BZ/2 \BZ \rightarrow \BZ/4\BZ.$ 
Let $\rad$ be the radical  of $\cdot.$ We say that $q$ is proper if $q$ vanishes on $\rad.$ If $q$ is proper,
the Brown invariant $\beta(q) \in \BZ/8\BZ$ is defined by the equation:

\[ 
 e^{\frac {2 \pi i \beta(q)} 8} = \frac 1 {\sqrt{2}^{\dim (V) + \dim (\rad) }} \sum_{v \in V} i^{q(v)} .
\]
We remark that the Brown invariant is additive for the direct sum of quadratic refinements. A simple graphical scheme for writing a given  form as a direct sum of simple elementary forms and thus calculating the Brown invariant is given in \cite{G2,G3}.

Let $L$ be a link in  $S^3$ with oriented components $\{K_i\}$.
Let $F$ be a not necessarily orientable spanning surface for $L$. 
Define a map
\[ q_F:H_1(F, \BZ/2 \BZ) \rightarrow \BZ/4\BZ\]
as follows. Given $x \in H_1(F,\BZ/2 \BZ)$, pick a simple closed curve  $\alpha_x$
representing $x$, and define $q_F(x)$  to be  the number of positive half twists in  a tubular neighborhood of $x$ in  $F$. More precisely,
\[ q_F(x)= \lk(\alpha_x, \hat {\alpha_x}) .\]
Here  $ \hat {\alpha_x}$ is the boundary of a tubular neighborhood of $\alpha_x$ oriented 
in the same direction as some arbitrarily chosen orientation for $x$. $q_F$ is well defined and is a quadratic refinement   
of the intersection pairing:
\[ \cdot : H_1(F, \BZ/2 \BZ) \times H_1(F, \BZ/2 \BZ)  \rightarrow \BZ/2\BZ .\]

We note that 
\[ q_F( [K_i] ) = 0 \text{ if and only if }  2 \lk(K_i, L-K_i) \equiv  0 \pmod{2}. \]
Let $\bar F$ denote the quotient space of $F$ obtained by identifying each component of $F$ to a point corresponding to that component. Let $\pi:  F \rightarrow \bar F$ denote the quotient map. Thus  $L$ is proper if and only if $q_F[K_i]=0 $  for all  $i$ if and only if $q_F$ can be factored through $\pi_*.$ 

Define
\[ \mu(F) = \frac 1 2 \sum_{i,j} \lk (K_i, K'_j) \]
where $K'_j$ denotes $ K_j$ pushed slightly into the interior of $F$. Then define
\[  \Arf (L) = \beta(q_F) - \mu(F) \in \BZ/ 8\BZ .\]

 One can relate any two spanning surfaces by a sequence of  moves  of certain types (and their inverses):    0)  isotopy  1) add a hollow handle and 2) replace a collar of the boundary with a punctured Mobius band, 3) the birth an empty two-sphere in a small 3-ball,  as described in \cite{G1}. The equivalence relation  generated by these moves is called $S^*$-equivalence.
Since $\beta(q_F) - \mu(F) $ is preserved by these moves,  $\Arf(L)$ is well defined. It takes values in $4\BZ/ 8\BZ$. It is invariant under oriented band summing, and adding and removing small unlinked unknots. This implies that it is an invariant of planar cobordism \cite{G2}.

\begin{figure}[ht]
\includegraphics[width=2in]{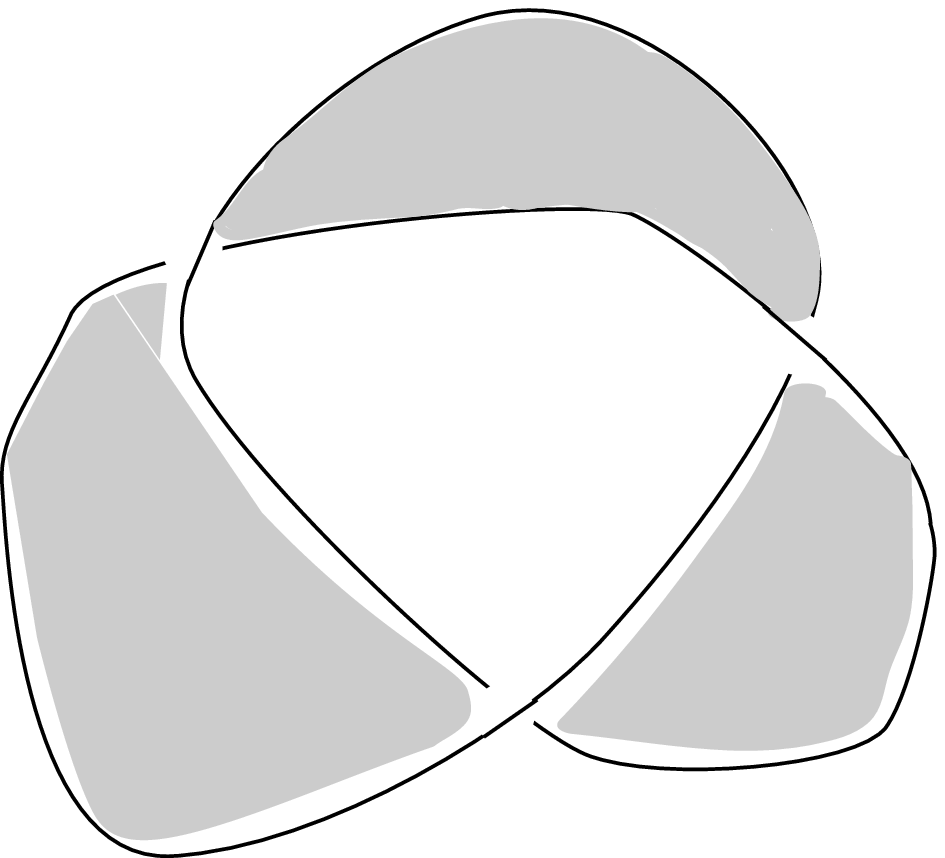}
\label{trefoil}
\end{figure}

Consider the trefoil in Figure \ref{trefoil} which is spanned by a Mobius band $F$. 
$H_1(F,\BZ/2 \BZ) \approx \BZ/2 \BZ$ and is generated by the core $x$.  We have
$q_F(x) \equiv -3 \equiv 1 \pmod{4}. $ $ \beta(q_F) \equiv  1 \mod{8},$
$\mu = \frac {-6} 2 =-3$, $\Arf( \text{trefoil}) \equiv 4 \pmod{8}.$

\section{Arf invariants of links in rational homology spheres}

Let $L$ be a link in rational homology sphere $M$ with oriented components $\{K_i\}$. 
By a spanning surface for  $L$, we mean a possibly non-orientable surface $F$ in $M$ with boundary $L$.  We may speak of $S^*$-equivalence of spanning surfaces in $M$.
Two different spanning surfaces for a link  need not be $S^*$-equivalent.  We introduced a parameter to index $S^*$ equivalence classes of spanning surfaces in \cite{G2}.  Suppose the  homology class of $L$ represents zero in $H_1(M,\BZ/ 2 \BZ).$  Then 
\[ \Gamma(L) = \{ \gamma \in H_1(M) | 2 \gamma = [L] \in H_1(M) \}\]
is nonempty. 

Let $F$ be a spanning surface for $L$ with no closed components, and $i_F:F \ra M$ be the inclusion. We have that $H_1(F)$ is free abelian and the homology class of $L$ the boundary of $F$ equipped with the string orientation of the link represents a homology class in $H_1(F)$ which is divisible by two. Define 

\[ \gamma(F) = i_F( \frac 1 2 [L] \in H_1(F) ) \in  \Gamma(L) \subset H_1(M) . \]

If $F$ has some  closed components, we take $ \gamma(F)$ to be $ \gamma(F')$ where $F'$  is formed by deleting the interior of a  disk from each closed component and orienting the new boundary components in an arbitrary way. For further discussion, see
  \cite{G1}.
  
We have that  $\gamma(F)$ is preserved by $S^*$-equivalence, and the map from $S^*$-equivalence classes of spanning surfaces for $L$ to $\Gamma(L)$ is bijective. By abuse of notation, we let $F_\gamma$ denote a spanning surface $F$ with $\gamma(F_\gamma)=\gamma$. This should not cause any confusion.

The number of positive half twists in a tubular neighborhood of a curve on a surface in $M$  has no well-defined analog in this more general situation.  What we actually need is an analog of the   number of half twists modulo  four.  

We do have a linking number in $\BQ$ which can be defined for disjoint 1-cycles. We continue to denote this linking number by $\lk$. We do know what it means to increase or decrease the number of half twists in the neighborhood of a curve $\alpha$.  Moreover if we add a half twist to a neighborhood of $\alpha$, we increase the linking number of  $\alpha$ and the boundary of the neighborhood of $\alpha$ by one. 

But we do not know what untwisted (mod 4) should mean.
As a replacement for this, it suffices to fix a quadratic refinement of the linking form $M$ \cite{G3}. The linking form of $M$,
\[ \ell_M: H_1(M) \times H_1(M) \ra \BQ/\BZ ,\] is a bilinear form with an injective adjoint.
A quadratic refinement of $ \ell_M$ is a function 
\[ r: H_1(M) \ra \BQ/\BZ ,\] such that
\[r(x+y) -r(x) -r(y) \equiv \ell_M(x,y) \pmod{1} .\]

It follows that the boundary $\hat \alpha$  of a neighborhood of an oriented  curve  $\alpha$ on a surface $F$ will have
$\lk(\hat \alpha, \alpha ) \equiv  r(3[\alpha])- r(2[\alpha]) -r([\alpha]) \equiv (3^3-2^2 -1) r([\alpha]) \equiv 4 r([\alpha])\pmod{1}.$  Thus we may define:
\[ q_{r,F}: H_1(F,\BZ/2\BZ) \rightarrow  \BZ/4 \BZ \]
by
\[q_{r,F}(x) \equiv \lk ( \gamma_x, \hat  \gamma_x ) - 4 r(x) \pmod{4} . \]
This is a well defined \cite[Thm. 6.1]{G2} quadratic refinement of the intersection form on $H_1(F,\BZ/2\BZ)$.

We also have that \cite[Prop 6.3]{G2}
\[ q_{r,F}( [K_i] ) = 0 \text{ if and only if } \lk(K_i, L-K_i) \equiv  2  \ell_M( [K_i],\gamma(F)) - 2 r([K_i]) \pmod{2} .\]
Note that $q_{r,F_\gamma}[K_i]=0 $ for all $i$ if and only if $q_{r,F_\gamma}$ can be factored through $\pi_*.$ 
 We say $(L,\gamma, r) $ is  { \em proper} if either of these equivalent conditions holds, i.e. if  \[\lk(K_i, L-K_i) \equiv   2 \ell_M( [K_i],\gamma(F)) - 2 r([K_i]) \pmod{2}  \text{ holds for all $i$}.\]

As in the previous section, we define
\[ \mu(F_\gamma) = \frac 1 2 \sum_{i,j} \lk (K_i, K'_j) \]
where $K'_j$ denotes $ K_j$ pushed slightly into the interior of $F_\gamma$.

If $(L,\gamma, r) $ is proper,  we define:
\begin{equation}  \Arf (L,\gamma,r ) = \beta(q_{r,F_\gamma}) - \mu(F_\gamma) \in \BQ/ 8\BZ .\label{Adef} \end{equation}

As for proper links  in $S^3$, $\Arf$ is well defined and is an invariant of planar cobordism.
 It is also invariant under oriented band summing, and adding and removing small unlinked unknots. 
When defined,   $ \Arf (L,\gamma,r )$, taken modulo four, depends only on  $\gamma$ and $r$
\cite[Prop. 6.9]{G2}.

 We remark that in \cite{G2,G4} our $\mu$ above is broken  up into the sum of two terms:
\[ \lambda(L) =  \sum_{i<j} \lk (K_i, K_j) \] and 
\[ \frac {-1} 2 e(F_\gamma)= \frac 1 2 \sum_{i} \lk (K_i, K'_i) .\]
For our purposes in this paper,  the use of $\mu$ makes things simpler.
We note that 
\[ \mu(F_\gamma) \equiv \ell([L], \gamma) \pmod{1}. \]
This restricts the range of values $ \Arf (L,\gamma,r )$ can take.

Both $\Gamma(L)$ and the set of all quadratic refinements of $\ell(M)$
are free $H^1(M,\BZ_2)$  sets. Properness of links  is preserved under the action which changes both $\gamma$ and $r$ simultaneously by  \cite[6.5]{G4}.  Moreover  $\Arf (L,\gamma,r)$ changes in a nice way under this action \cite[6.5]{G4}:

\begin{prop} \label{change}  If  $( L,\gamma,  q)$ is proper then $(L,\psi \cdot \gamma,   \psi \cdot q)$ is proper and \[\Arf (L,\psi \cdot \gamma,   \psi \cdot q) = \Arf ( L,\gamma,  q) +\Arf( \emptyset, \psi \cdot 0,   \psi \cdot q) . \] \end{prop}

\section{the tangent circle bundle  of the real projective plane}

We will be concerned with links in the tangent circle bundle of
$\rp$ which we will denote $\T$. This is the boundary of $\D$, the tangent disk bundle of $\rp$. $\D$, of course,  has the homology of $\rp$. It follows that $\T$ is a rational homology sphere. In fact
$\T$ is the lens space $L(4,3)$ \cite{G1}. We let $\ell$ denote $\ell_{\T}$. 

By an  $\I$-curve, we mean a collections of immersed oriented curves in $\rp$ where each curve is never tangent to itself or any other curve in the collection. 
We can describe links in $\T$ by $\I$-curves. 
Lying over an $\I$-curve $C$ in $\T$, we have a link $\T(C)$ consisting of all the rays tangent to a part of $C$ and pointing in the direction of the orientation.  We let $g$ denote the homology class of $\T$ represented by the  oriented knot in $\T$ 
which lies over a straight line with either orientation.\footnote{ There is an isotopy that can be obtained by  spinning the line around a point on the line. This isotopy  sends a line in $\rp$ with  one orientation to the same line with the opposite orientation.} The homology class, $g$ is a generator for $H_1(\T)$.  An oval represents $2g$.

Whenever we have an invariant of links in $\T$, we obtain an invariant of  $\I$-curves.
We use the same symbol to denote an invariant of oriented links in $\T$ and the corresponding invariant of  $\I$-curves.  An immersed circle in an $\I$-curve
$C$  is called a component of $C$, and describes a component of the link  $\T(C).$ This should not cause  any  confusion and simplifies our 
expressions.

Similarly, if we have two $\I$-curves $C_1$ and $C_2$ which are never tangent to  each other, we can speak of their linking number $\lk(C_1, C_2) \in \frac 1 4 \BZ.$
In  \cite{G1}, we  worked out the linking numbers of some two component $\I$-curves.  
Two linked ovals (see \S1) which are oriented in the same  direction (so they contribute to $\Pi^-$) have linking number one. Two linked ovals (see \S1) which are oriented oppositely (so they contribute to $\Pi^+$) have linking number minus one. Two unlinked ovals, have linking number zero.

 A one-sided simple closed curve and a disjoint oval have linking number $\frac 1 2$, if the oval  is homologous to twice the one-sided curve in the Mobius band formed by deleting the interior of the oval from $\rp$. Otherwise the linking number is $-\frac 1 2$.  A one-sided simple closed curve and  oval which the curve meets twice transversally have linking number $-\frac 1 2$. Two one-sided simple closed curves which meet in one point have linking number $-\frac 14$.  Thus $\ell(g,g) = - \frac 1 4,$ and there are two quadratic refinements of $\ell$, $r_{\frac {-1} 8}$ with $r_{\frac{ -1} 8}(g)= 
  \frac{ -1}  8$, and, $r_{ \frac 3 8}$ with $r_{ \frac 3 8}(g)=   \frac 3 8$. 

A  dangerous type two move \footnote{ A dangerous type two move is a non-safe type two Reidemeister move on diagrams}  between two components which reduces the 
number of double points by two leads to a new  $\I$-curve where the linking number between the two components has been increased by one.  Smoothing (according to the orientations) a double point of an $\I$-curve is called a smoothing move.  A smoothing move  corresponds to a oriented band move to the corresponding link in $\T$.  Thus  if we wish to calculate a linking number between two sub-$\I$-curves of an $\I$-curve, we can  smooth the double points of sub-$\I$-curves  without changing the relevant linking numbers. Of course the reverse of this smoothing move which we call an unsmoothing  move  will introduce a double point and also corresponds to  an oriented band move performed on the corresponding link in $\T$.

Lying above an empty figure eight curve is a local unknot in $\T$. The easiest way to see this is to note  that any empty figure eight  is isotopic to an  $\I$-curve which does not use every tangent direction as one travels around a circuit.  Thus this  $\I$-curve lies in an 
$I$-bundle over a disk  in $\rp$ which includes the figure eight. Moreover the $\I$-curve projects to a curve  in this  disk with only  a one double point.
Also we have that  lying above two $\I$-curves  related by balanced type I move,  a safe type $II$ or a type $III$ move are isotopic links in $\T$. Thus weakly equivalent $\I$-curves are proper for that same $\gamma$ and $r$, and have the same $\Arf$ invariant, when it is defined.

As an example, consider the $\I$-curve, called $\td$ given by the  diagram for trefoil in Figure 2, with the three crossings (made double points) placed in an affine part of $\rp$.  A Mobius band $F$ is described by a ``vector field'' \footnote{At the double points, we have a whole arc of lines joining the  two intersecting lines which fill out the shaded region.} on the shaded region which extends the tangential field on the boundary. Let $\delta$  denote a 1-sided curve in $F$ which  represents the homology class $2g.$  Thus  $\gamma(F) = 2g$. The linking number of the boundary of $F$ with a parallel may be calculated using the above techniques. It is $-2$. Thus $\mu(F) =-1.$ . Moreover it follows that 
\be q_{ \frac {-1} 8, F}(\delta) = q_{ \frac {3} 8, F}(\delta)=1 . \label{delta} \ee  Thus we have 
$\Arf(\td,2g, \frac {-1} 8)\equiv \Arf(\td,2g, \frac {3} 8)\ \equiv 2 \pmod{8}$. $\td$ is actually a nodal M-curve of degree 4. Thus this is a small confirmation of Theorem \ref{main} below.

\begin{prop}\label{changeprop}   \[\Arf (\emptyset, 2g, r_{  \frac {-1} 8}) \equiv  2 \pmod{8} \text{ and} \]
\[\Arf (\emptyset, 2g, r_{ \frac 3 8}) \equiv  -2 \pmod{8} . \]
If  either $(L, \gamma,r_{  \frac {-1} 8})$ or $(L, \gamma +2g ,r_{ \frac 3 8})$ is proper for some $\gamma$,  then the other is proper, and
\[ \Arf(L, \gamma,r_{ -\frac {1} 8})-\Arf(L, \gamma + 2g ,r_{ \frac 3 8}) \equiv   2 \pmod{8}.\]
\end{prop}

\begin{proof} Fix a line  in $\rp$ and consider the Klein bottle $K$ in  $\T$ given by the set of all directions through points on this line. $K$ is a spanning surface for the  empty link. The lift of this line  with one orientation $\alpha$, and the lift $\alpha'$  with the other orientation each have neighborhoods which are Mobius bands. These neighborhoods may be isotoped off of  $K$ so that they are the two lifts to $\T$ of the Mobius bands neighborhood of the line  in   $\rp$.
Since $\lk( \alpha, \hat \alpha)$  and  $\lk( \alpha', \hat \alpha')$ are seen to be  $\frac 1 2$,  we  calculate that: 
 \be q_{r_{ \frac {-1} 8},K}( \alpha)=  q_{r_{ \frac {-1} 8},K}( \alpha') = \frac 12 -4 \cdot{ \frac {-1} 8} = 1 \text{ and}  \label{alpha1} \ee
\be  q_{r_{ \frac 3 8},K}( \alpha)=  q_{r_{ \frac 3 8},K}( \alpha')=  \frac 12 -4 \cdot { \frac {3} 8} =  -1  .
 \label{alpha2}
\ee
 In this way we obtain the first two equations. The rest follows from Proposition \ref{change}.
\end{proof}

 The following is a special case of \cite[Prop 3.2]{G5}. Its proof is a simple calculation using as a spanning surface the union of $k$ disjoint annuli swept out by rotating $k$  lines.

\begin{prop}\label{X} Let $X_{2k}$ be $2k$ lines straight lines in general position, then  
$(X_{2k},kg, r_{  \frac {-1} 8})$ and $(X_{2k},(k+2)g, r_{  \frac {3} 8})$ are  proper, and  
\[ 
\Arf (X_{2k}, kg, r_{  \frac {-1} 8}) 
\equiv \frac {k^2} 2 
\equiv 
\Arf (X_{2k}, (k+2) g, r_{\frac {3} 8}) +2  \pmod{8}. \]
Also $(X_{2k},(k+2)g, r_{  \frac {-1} 8})$ and $(X_{2k},kg, r_{  \frac {3} 8})$ are  not proper.
\end{prop}

If $C$ is a collection of disjoint ovals in $\rp$, let $B^+(C)$
be the closed surface in $\rp$ with boundary $C$ which is the closure of the set of points which lie  inside an odd number of ovals of $C.$
Let $B^-(C)$ be  the  closure of  $\rp \setminus B^+(C)$.

\begin{prop} \label{odd} Let $C$ be a collection of disjoint ovals in $\rp$ with an even number of components,
 then  the following are equivalent;
 \begin{enumerate}
 \item $(C,\gamma,r) $ is proper for some $\gamma$, and $r$.
 \item $(C,\gamma,r) $ is proper for  all possible $\gamma$ (i.e.  both $0$ and $2g$ ) and for all $r$.
 \item $C$ is odd.
  \item Every component of $B^+(C)$ has even Euler characteristic.
 \end{enumerate}
 \end{prop}
 \begin{proof}  Each component $C_i$ represents $2g\in  H_1(\T)$.
 Thus $2r([C_i])=8r[g]=1 \pmod{2}.$
 As $C$  has an even number of components and each oval represents $2g$, we have that $[C]= 0 \in  H_1(\T)$. Thus $\gamma$ must be either $0$ or $2g$. Also $C$ is proper for either $\gamma$ and either $r$ if and only  if $\lk(C_i,C-C_i)$ is odd. On the other hand,
 $\lk(C_i,C \setminus C_i)$ is odd if and only if $C$ is odd. The equivalence of the last two conditions is easily seen.
\end{proof}

The  definitions of this paragraph are due to Rokhlin \cite{R}. An odd oval is called disoriented if forms a  negative pair with  the  even oval that immediately surrounds it. Let $d$ denote the number of disoriented ovals. Let $D^+$ denote the number of positive pairs with disoriented outer oval. Similarly let $D^-$ denote the number of negative pairs with disoriented outer oval. Rokhlin observed that
\[\Pi^+ -\Pi^- = n -2( d-D^+ + D^-) .\]
If $C$ is odd, it is easy to see that $D^+ + D^-$ is even. Thus,  if $C$ is odd, we have that
\begin{equation} 2d \equiv \Pi^+ -\Pi^- -n \pmod{4} . \label{2d}\end{equation}

\begin{prop}\label{oddf} Let $C$ be an odd collection of disjoint ovals in $\rp$, then  
$(C, (\Pi^+ -\Pi^- -p )g, r)$ and $(C, (\Pi^+ -\Pi^- -p -2 )g, r)$ are proper for either $r$.   Moreover 
\[\Arf (C, (\Pi^+ -\Pi^- -p )g, r_{ \frac {-1} 8} )\equiv \Arf (C, (\Pi^+ -\Pi^- -p )g, r_{ \frac {3} 8} )\equiv \Pi^+ -\Pi^- -p \pmod{8} , \]
\[\Arf (C, (\Pi^+ -\Pi^- -p-2 )g, r_{\frac 3 8)}\equiv \Pi^+ -\Pi^- -p-2 \pmod{8}  \text{ and} \]
\[\Arf (C, (\Pi^+ -\Pi^- -p+2 )g, r_{ \frac {-1} 8)}\equiv \Pi^+ -\Pi^- -p+2 \pmod{8} . \]
\end{prop}

\begin{proof} By Proposition \ref{odd}, every component of $B^+$ has even Euler characteristic. Thus we can pick a vector field on $B^+$ which is tangent to the boundary and pointed in the direction of the orientation of $C$ with $(n-p)/2$ zeros of index -2. This defines a surface in $\T$ lying over $B^+$ with $(n-p)/2$  disks removed around the zeros, which we denote by $z_i$.  We complete this surface to form a  spanning surface  $F$ for $L(C)$ in  $\T$ by adding Mobius bands above the disks around the zeros. We do this  so that  the cores of  these Mobius bands, which we denote by $x_i$, are the fibers over $z_i$ of the map from $F$ to $\rp$. This is similar to the construction of spanning surfaces  in  \cite[\S1]{G4}.  Then  
\be q_{r,F}( x_i  ) \equiv -1 -4  r(2g)\equiv 1 \pmod{4} . \label{x} \ee
Since  $B^+$ is a planar surface, \[\beta(q_{r,F}) \equiv \frac {n-p} 2 \pmod{8}.\]
We have that
\be \label{mu}  \mu(F) = \frac 1 2( n+p-2(\Pi^+ - \Pi^-)). \ee
Let $d$ be the number of even ovals  which form a  negative pair with  the  oval that immediately surrounds them.
We have  that 
\[ \gamma(F) = (\frac {n-p} 2 +d )\  2g.\]
Using Equation \ref{2d}:
\[ \gamma(F) =  ( \Pi^+ - \Pi^- -p)\ g.\]
Together with Equation \ref{Adef} gives the first equation. The last two equations  follow from this and 
Proposition \ref{changeprop}
\end{proof}

 We have:
 
 \begin{prop} \label{even} Let $C$ be a collection of disjoint ovals in $\rp$ with an odd number of components,
 then  the following are equivalent;
 \begin{enumerate}
 \item $(C,\gamma,r) $ is proper for some $\gamma$, and $r$.
 \item $(C,\gamma,r) $ is proper for all possible $\gamma$ (i.e.  both $g$ and $-g$ ) and for all $r$.
 \item $C$ is even.
  \item Every component of $B^-(C)$ has even Euler characteristic.
 \end{enumerate} 
\end{prop}
\begin{proof} 
 As $C$  has an odd number of components and each oval represents $2g$, we have that $[C]= 2g \in  H_1(\T)$. Thus $\gamma$ must be either $g$ or $-g$. Also $C$ is proper for either $\gamma$ and either $r$ if and only  if $\lk(C_i,C-C_i)$ is even. On the other hand,
 $\lk(C_i,C \setminus C_i)$ is even if and only if $C$ is even. The equivalence of the last two conditions is easily seen.
\end{proof}

\begin{prop}\label{evenf} Let $C$ be an even collection of disjoint ovals in $\rp$, then  
$(C, \pm g, r)$  are proper for both $r=r_{ \frac {-1} 8}$ and $r=r_{ \frac {3} 8}$.   Moreover 
\[\Arf (C, \pm g,r_{ \frac {-1} 8})\equiv \Pi^+ -\Pi^- -n +\frac 1 2 \pmod{8} , \]
\[\Arf (C, \pm g, r_{ \frac {3} 8})\equiv \Pi^+ -\Pi^- -n - \frac 3 2 \pmod{8} . \]
\end{prop}

\begin{proof} By Proposition \ref{even}, every component of $B^-$ has even Euler characteristic. Thus we can pick a vector field on $B^-$ which is tangent to the boundary and pointed in the direction of the orientation of $C$ with $(p-n-1)/2$ zeros of index -2. This defines a surface in $\T$ lying over $B^-$ with $(p-n-1)/2$  disks removed around the zeros.  As above, we complete this surface to form a  spanning surface  $F$ for $L(C)$ in  $\T$ by adding Mobius bands above the removed disks around the zeros.
$B^-$ is a planar surface disjoint union a Mobius band with some holes removed. The same can be said of the surface we obtain when we delete neighborhoods of the singularities of the vector field.
 One can see the $\gamma(F)$ is either $g$ or $-g$. As we will see,  this allows us to compute our Arf invariants without ambiguity. Let us now  write $\gamma(F)= \pm g$, and read plus, if indeed it is plus, and read minus, if indeed it is minus. 
 
As in proof of  Proposition \ref{oddf}, we  let $x_i$ denote  the cores of  these Mobius bands.
and Equation \ref{x} gives us $q_{r,F}(x_i )$.
Let $\alpha$ denote an orientation reversing curve in  $B^-$, then, as in Equations \ref{alpha1} and \ref{alpha2}, 

\begin{equation}  \label{alpha} q_{r,F}(\alpha) \equiv \frac 1 2  -4  r(g)\equiv 
\begin{cases} 
+1  , &\text{ if $r =r_{\frac {-1} 8}$}\\
-1 . &\text{ if $r =r_{\frac {3} 8}$} 
\end{cases} \pmod{4}. 
 \end{equation}

Thus  \[\beta(q_{r,F}) \equiv  \frac {p-n-1} 2  +  \begin{cases} 
+1  , &\text{ if $r =r_{\frac {-1} 8}$}\\
-1 . &\text{ if $r =r_{\frac {3} 8}$}
\end{cases}   \pmod{8}.\]
Equation \ref{mu} gives $ \mu(F)$.
Plugging into the  equation \ref{Adef},  we obtain the stated results but where must read $\pm$ according to whether
$\gamma(F)$  is $\pm g$. However an application of  Proposition \ref{changeprop}
 to both these equations shows that they must hold for the other choice of $\gamma$ as well.
 \end{proof}

\section{Main results}

By the proof of  \cite[Theorem 3.1]{G4}  :  

\begin{thm}\label{cob}Let $\A$ be a nodal M-curve of degree $2k$,  then  $\T(\A)$ is planar cobordant to $\T(X_{2k}).$
\end{thm}

We have the following simple corollary. 

\begin{prop}\label{simple.cor} Let $\A$ be a nodal M-curve of degree $2k$.   $\T(\A)$ is homologous to $2kg \in H_1(\T).$ If $k$ is  even,
$\A$ cannot be weakly equivalent to an even simple curve.  
If $k$ is  odd, $\A$ cannot be weakly equivalent to an odd simple curve.
\end{prop}

The next Theorem follows from Proposition \ref{X}, the invariance of Arf invariants of proper links under planar cobordism, and Theorem
\ref{cob}.

\begin{thm}\label{main} Let $\A$ be a nodal M-curve of degree $2k$, and  suppose
$(\A,kg,  r_{\frac {-1} 8})$ is proper, then  $\Arf (\A, kg,  \ r_{\frac {-1} 8}) \equiv \frac {k^2} 2 \pmod{8}$.
\end{thm}

We can now give the proof of Theorem \ref{simple.thm}. Theorem \ref{main} is more general but its application requires that one calculate $\Arf (\A, kg,  r_{\frac {-1} 8})$. 

\begin{proof}[Proof of Theorem \ref{simple.thm}]  Suppose $k$ is even, and $\A$ is weakly equivalent to an odd simple curve.  Then, by Proposition \ref{odd}, $(\A,kg,   r_{\frac {-1} 8})$ is proper. 
So, by Theorem \ref{main}, $\Arf (\A, 0,  r_{ \frac {-1} 8}) \equiv  \frac {k^2} 2 \pmod{8}$. By Proposition \ref{oddf}, one of $(\Pi^+ -\Pi^- -p )g$ or $(\Pi^+ -\Pi^- -p +2)g$ is zero, and
\[\Arf (C, 0, r_{ \frac {-1} 8} )\equiv  \Pi^+ -\Pi^- -p + \begin{cases} 
0  \pmod{8}  , &\text{ if $\Pi^+ -\Pi^- -p \equiv 0 \pmod{2} $}\\
2   \pmod{8}  &\text{ if $\Pi^+ -\Pi^- -p  \equiv 2 \pmod{2} .$}
\end{cases} \]
Since $\T(C)$ and $\T(\A)$ are planar cobordant, we have 
\[\Arf (\A, 0,  r_{ \frac {-1} 8}) \equiv \Arf (C, 0, r_{ \frac {-1} 8} )\pmod{8}.\]
This gives the $k$ even case.

Now suppose $k$ is odd, and $\A$ is weakly equivalent to an even simple curve.  Then, by Proposition \ref{even},  $(\A,kg,   r_{\frac {-1} 8})$ is proper. 
So, by Theorem \ref{main}, $\Arf (\A, kg,  r_{ \frac {-1} 8}) \equiv  \frac {k^2} 2 \pmod{8}$. By Proposition \ref{evenf}, 
\[\Arf (C, k g, r_{ \frac {-1} 8} )\equiv  \Pi^+ -\Pi^- -n+ \frac 1 2 \pmod{8} .\]
Since $\T(C)$ and $\T(\A)$ are planar cobordant, we have 
\[\Arf (\A, kg,  r_{ \frac {-1} 8}) \equiv \Arf (C, kg , r_{ \frac {-1} 8} )\pmod{8}.\]
This gives the $k$ odd case.

\end{proof}

We now wish  to show how Fiedler's original congruence for certain nonsingular curves is equivalent to Theorem \ref{simple.thm} when
 $\A$ is a simple curve, even or  odd, as the case may be.

 \begin{thm}[Fiedler]\label{fied} Let $\A$ be a nonsingular  M-curve of degree $2k$. 

If $k$ is even, and $\A$ is  an odd simple curve, then
\[ 
p -n \equiv -k^2
\pmod{16} .\]

If $k$ is odd, and $\A$ is an even simple curve, then
\[ 
p-n \equiv  
 1
 \pmod{16}. \]
\end{thm}

\begin{proof} We have Gudkov's congruence:
\be p-n \equiv k^2 \pmod{8}.
\label{gu}
\ee

Harnack's inequality is extremal:

\be
p+n = 1 + \binom{2k-1}{2} =2 k^2-3k +2 .
\label{ha}
\ee

Adding  Equations \ref{gu} and \ref{ha} , and dividing by $2$:
\be
p \equiv \frac {3 k^2-3k+2}{2} \pmod{4} .
\label{guha}
\ee

According to Rokhlin \cite[Equation 4]{R}, 
\be
\Pi^+ - \Pi^- = \frac {(k-1)(k-2)}{2}= \frac{k^2-3k+2}{2}.
\label{ro}
\ee

Subtracting Equation \ref{guha} from Equation \ref{ro}, we obtain
\be
\Pi^+ - \Pi^- -p \equiv -k^2 \pmod{4} .
\label{init}
\ee 

At this point, we consider separately three different cases:
$k \equiv 0 \pmod{4}$, $k \equiv 2 \pmod{4}$, and $k \equiv 1 \pmod{2}$. 

Assume now that $k \equiv 0 \pmod{4}$, and $\A$ is an odd simple curve.  By Equation \ref{init}, 
$
\Pi^+ - \Pi^- -p \equiv 0 \pmod{4}.
$
Thus by Theorem \ref{simple.thm} 
\[
\Pi^+ - \Pi^- -p \equiv \frac {k^2} 2 \equiv 0 \pmod{8}.
\]
Thus using Equation \ref{ro},
\[ 2p = k^2-3k +2 \equiv -3k +2 \pmod{16}. \] 
Thus, subtracting  Equation \ref{ha},
\[ p-n = -2k^2 \equiv 0 \pmod{16} . \]
This agrees with the conclusion of the Theorem to be proved.

Now assume  that $k \equiv 2 \pmod{4}$, and $\A$ is an odd simple curve.  We still have
$
\Pi^+ - \Pi^- -p \equiv 0 \pmod{4}
$, but now $\frac {k^2} 2 \equiv 2 \pmod{8}$. So  
by Theorem \ref{simple.thm} 
\[
\Pi^+ - \Pi^- -p \equiv \frac {k^2} 2-2 \equiv 0 \pmod{8}.
\]
Using Equations \ref{ro} and \ref{ha}, as above, we obtain the conclusion of the Theorem to be proved.

Finally assume  that $k \equiv 1 \pmod{2}$, and $\A$ is an even simple curve. By  Theorem \ref{simple.thm} 
\[
\Pi^+ - \Pi^- -n \equiv \frac {k^2-1} 2 \equiv 0 \pmod{8}.
\]
Using Equation \ref{ro}:
\[ -2n \equiv 3k-3 \pmod{16}.\] Using Equation \ref{ha}, we obtain
\[ p-n \equiv  2 k^2 -1 \equiv 1 \pmod{16} .\]

\end{proof}

\section{Fiedler's Curves}

\subsection{Prohibiting a curve of degree 6} 
\begin{figure}[ht]
\includegraphics[width=1.5in]{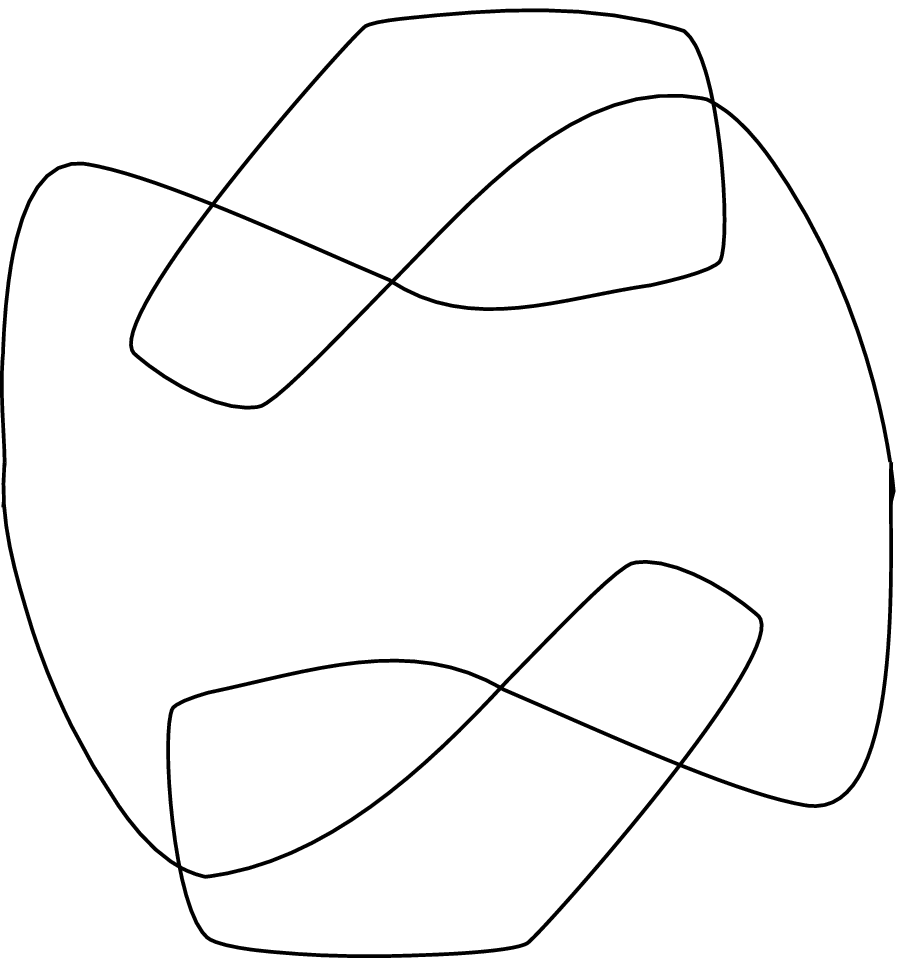}
\caption{$C_2$} \label{deg6}
\end{figure}

 Consider the hypothetical curve of degree six prohibited by  Fiedler \cite[Figure 1]{F2}. We denote this   $\I$-curve by $C_1$.  We note that there is only one orientation on this two component curve up isotopy. So we equip $C_1$ with the orientation which allows the unsmoothing move that we take below. We have that $C_1$ is proper for $\gamma = \pm g$ and $r$ equal either  $r_{\frac {-1} 8}$ or $r_{\frac {3} 8}.$  We perform  some safe type two moves and type three moves  on $C_1$  followed by an unsmoothing  move, and a balanced type one move to obtain the $\I$-curve in Figure \ref{deg6}, which we denote by $C_2$. Since the unsmoothing move decreases the number of components, $(C_2,\pm g, r)$ must be proper, \cite[Corollary 6.10]{G2}. So we have:
 \[ \Arf(C_1,\pm g, r) \equiv  \Arf(C_2,\pm g, r) \pmod{8} \]
 for either $r$.
 
  We pick a ``vector  field'' on $B^-(C_2)$ which is tangential to the boundary  pointed in the direction of the  orientation, has   a whole arc of tangent directions at each double point (as in the spanning surface for $\td$ in section 3 ) and has a single singularity of index  $-2$. This describes a spanning surface $F$ for $L(C_2)$. We calculate that $\gamma(F)= -g$ and 
 $\mu(F) = \frac {-9} 2.$  We  have a basis for $H_1(\bar {F}, \BZ/2 \BZ)$ consisting of two 1-sided curves $\delta_1$, 
 $\delta_2$ on $F$ in the affine part of picture, the  fiber $x$ over the singularity of index $-2$,  and the line at infinity $\alpha$. As in Equation \ref{delta}
 \[ q_{ \frac {-1} 8, F}(\delta_1)= q_{ r_{\frac {-1} 8}, F}(\delta_2)=1 .\]
 As in Equation \ref{x}, \[q_{r_{ \frac  {-1} 8},K}(x) = 1.\]
As in Equation \ref{alpha},
\[q_{r_{ \frac  {-1} 8},K}( \alpha) = 1.\]
Thus \[\beta( q_{r_{ \frac  {-1} 8},F}) \equiv 4 \pmod{8} .\]
So \[\Arf(C_2, 3g, r_{ \frac  {-1} 8}) \equiv \Arf(C_2,-g, r_{ \frac  {-1} 8}) \equiv 4+ \frac 9 2 \equiv  \frac 1 2 \nequiv  \frac 9 2 \pmod{8}.\]
 Thus by Theorem \ref{main}, $C_1$ is not a real algebraic curve of degree 6.

\subsection{Prohibiting a curve of degree 8} 

\begin{figure}[ht]
\includegraphics[width=2.5in]{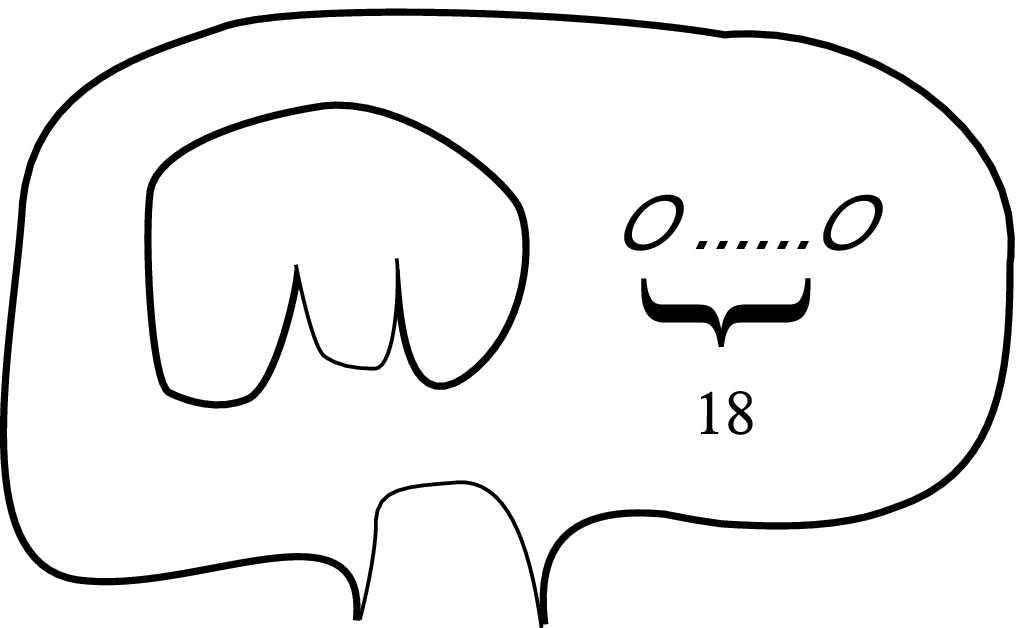}
\caption{The outer curve with the two flops is oriented clockwise. The inner curve with the two flops is oriented counterclockwise. Ten of the ovals are oriented counterclockwise, and eight of them are oriented clockwise.} \label{deg8}
\end{figure}

 Consider the hypothetical curve of degree eight prohibited by  Fiedler \cite[Figure 2]{F2}. We denote this   $\I$-curve by $C_3$ and equip it with the only semi-orientation on $C_3$, up to isotopy,  which is consistent with Rokhlin \cite[Equation 4]{R} and Fiedler's alternation of ovals with respect to a pencil of lines \cite{F3} when applied to $C_3$ smoothed. For this calculation, we use the language of  floppy curves as developed in  \cite{G4}.  We  draw a real floppy  curve $C_4$ in $\rp$ whose corresponding link in $\T$ is isotopic to the link that corresponds to $C_3$ with this orientation in Figure \ref{deg8}.   We check that 
 that $C_3$ is proper for $\gamma = 0$  and $r= r_{\frac {-1} 8}$.
 We can extend the vector field on the boundary over $B^+$ with  ten singularities of  index $-2$. This specifies  a spanning surface $F$. One calculates that  $\gamma(F) =0$ and $\mu(F)= 6$. Using equation \ref{x},  
  thus \[\beta( q_{r_{ \frac {- 1} 8},F})  \equiv 10 \equiv 2 \pmod{8} .\]
So \[\Arf(C_3,0, r_{ \frac {-1} 8}) \equiv  2-6 \equiv  4 \nequiv  0 \pmod{8}.\]
Thus by Theorem \ref{main}, $C_3$ is not a real algebraic curve of degree 8.

\end{document}